\documentclass[twosided,reqno]{amsart}

\usepackage{amsmath}
\usepackage{amssymb}
\usepackage{amsthm}

\theoremstyle{theorem}
\newtheorem{theorem}{\scshape Theorem }[section]

\theoremstyle{definition}

\numberwithin{equation}{section}

\begin{document}

\title[Degenerate Laplace transform and degenerate gamma function]{Degenerate Laplace transform and degenerate gamma function}

\author{Taekyun Kim}
\address{Department of Mathematics, Kwangwoon University, Seoul 139-701, Republic
of Korea\\
Department of Mathematics, College of Science,
 Tianjin Polytechnic University,
 Tianjin 300160, China}
\email{tkkim@kw.ac.kr}

\author{Dae San Kim}
\address{Department of Mathematics, Sogang University, Seoul 121-742, Republic
of Korea}
\email{dskim@sogang.ac.kr}



\maketitle

\begin{abstract}
In this paper, we introduce the degenerate Laplace transform and degenerate gamma function and investigate some properties of the degenerate Laplace transform and degenerate gamma function. From our investigation, we derive some interesting formulas related to the degenerate Laplace transform and degenerate gamma function.
\end{abstract}

\section{Introduction}
As is well known, the {\it{gamma function}} is defined by
\begin{equation}\label{1}
\Gamma(s)=\int_0 ^{\infty} e^{-t}t^{s-1}dt,{\text{ where }}s \in {\mathbb{C}}{\text{ with }}Re(s)>0,{\text{(see \cite{01, 02, 11, 27})}}.
\end{equation}
From \eqref{1}, we note that
\begin{equation}\label{2}
\Gamma(s+1)=s\Gamma(s), {\text{ and }}\Gamma(n+1)=n!,~(n\in{\mathbb{N}}).
\end{equation}
The {\it{degenerate exponential function}} $e_{\lambda} ^t$ is a function of two variables $\lambda \in (0,~\infty)$ and $t \in \mathbb{R}$ defined as
\begin{equation}\label{3}
e_{\lambda} ^t =(1+\lambda t)^{\frac{1}{\lambda}}.
\end{equation}
Note that
\begin{equation*}
\begin{split}
\lim_{\lambda\rightarrow0+}e_{\lambda} ^t=&\lim_{\lambda\rightarrow 0+} (1+\lambda t)^{\frac{1}{\lambda}}=\sum_{n=0} ^{\infty} \frac{1}{n!}t^n\\
= &e^t.
\end{split}
\end{equation*}

Let $f$ be a function defined for $t \geq 0$. Then the integral
\begin{equation}\label{4}
\mathcal{L}(f(t))=\int_0 ^{\infty} e^{-st}f(t)dt,{\text{ (see \cite{05, 07, 19, 20, 26, 27})}},
\end{equation}
is said to be the {\it{Laplace transform of $f$}}, provided that the integral converges.

When the defining integral \eqref{4} converges, the result is a function of $s$. As a matter of notation, if we use the lowercase letter to denote the function being transformed then the corresponding uppercase letter will be used to denote its Laplace transform: for example,
\begin{equation*}
\mathcal{L}(f(t))=F(s),~\mathcal{L}(g(t))=G(s),~\mathcal{L}(y(t))=Y(s),{\text{ and }}\mathcal{L}(h(t))=H(s),\ldots.
\end{equation*}
It is easy to show that
\begin{equation}\label{5}
\mathcal{L}(1)=\frac{1}{s},~(s>0),~\mathcal{L}(t)=\frac{1}{s^2},~(s>0),~\mathcal{L}(e^{-3t})=\frac{1}{s+3},~(s>-3),
\end{equation}
and
\begin{equation}\label{6}
\mathcal{L}(\sin 2t)=\frac{2}{s^2+4},~(s>0),~\mathcal{L}(\cos 2t)=\frac{s}{s^2+4},~(s>0).
\end{equation}

The {\it{Euler formula}} is defined by
\begin{equation}\label{7}
e^{i \theta}=\cos \theta +i \sin \theta,{\text{ (see \cite{25, 27})}}.
\end{equation}
Thus, by \eqref{7}, we get
\begin{equation}\label{8}
\cos a \theta=\frac{e^{ia \theta}+e^{-ia \theta}}{2},~\sin a\theta=\frac{e^{ia \theta}-e^{-ia \theta}}{2i}.
\end{equation}

From \eqref{3}, we consider the {\it{degenerate Euler formula}} which is given by
\begin{equation}\label{9}
e_{\lambda} ^{it}=(1+\lambda t)^{\frac{i}{\lambda}}=\cos_{\lambda} (t)+i\sin_{\lambda} (t).
\end{equation}
Note that
\begin{equation}\label{10}
\lim_{\lambda\rightarrow0}e_{\lambda} ^{it}=\lim_{\lambda\rightarrow0}(1+\lambda t)^{\frac{i}{\lambda}}=e^{it}=\cos t+i\sin t.
\end{equation}
By \eqref{9} and \eqref{10}, we get
\begin{equation}\label{11}
\lim_{\lambda\rightarrow 0} \cos_{\lambda}(t)=\cos t,~\lim_{\lambda\rightarrow 0}\sin_{\lambda}(t)=\sin t.
\end{equation}
From \eqref{9}, we note that
\begin{equation}\label{12}
\cos_{\lambda}(t)=\frac{e_{\lambda} ^{it}+e_{\lambda} ^{-it}}{2},~\sin_{\lambda} (t)=\frac{e_{\lambda} ^{it}-e_{\lambda} ^{-it}}{2i},
\end{equation}
which are called the degenerate cosine and degenerate sine functions, respectively.

Recently, several authors have studied the degenerate special polynomials and numbers (see [1-27]).

For example, L. Carlitz had studied degenerate Bernoulli polynomials given by the generating function
\begin{equation}\label{13}
\frac{t}{e_{\lambda} ^t-1}e_{\lambda} ^{xt}=\frac{t}{(1+\lambda t)^{\frac{1}{\lambda}}-1}(1+\lambda t)^{\frac{x}{\lambda}}=\sum_{n=0} ^{\infty}\beta_{n,\lambda}(x)\frac{t^n}{n!},{\text{ (see \cite{03})}}.
\end{equation}
When $x=0$, $\beta_{n,\lambda}=\beta_{n,\lambda}(0)$ are called {\it{degenerate Bernoulli numbers}} (see \cite{03, 04}).

In this paper, we introduce the degenerate Laplace transform and degenerate gamma function and investigate some properties of the degenerate Laplace transform and degenerate gamma function. From our investigation, we derive some formulas related to the degenerate Laplace transform and degenerate gamma function.

\section{Degenerate gamma function}

For each $\lambda\in(0,~\infty)$, we define the {\it{degenerate gamma function}} for the complex variable $s$ with $0< Re(s) <\frac{1}{\lambda}$ as follows:
\begin{equation}\label{14}
\Gamma_{\lambda}(s)=\int_0 ^{\infty}(1+\lambda t)^{-\frac{1}{\lambda}}t^{s-1}dt=\int_0 ^{\infty}e_{\lambda} ^{-t}t^{s-1}dt.
\end{equation}
Here we note that
\begin{equation}\label{15}
\begin{split}
\int_0 ^{\infty}(1+\lambda t)^{-\frac{1}{\lambda}}t^{s-1}dt=& \lambda^{-s}\int_0^{\infty}(1+t)^{-\frac{1}{\lambda}}t^{s-1}dt\\
=&\lambda^{-s}B(s,~\frac{1}{\lambda}-s),
\end{split}
\end{equation}
where $B(x,~y)$ is the Beta function given by
\begin{equation}\label{16}
B(x,~y)=\int_0^{\infty}\frac{u^{x-1}}{(1+u)^{x+y}}du=\frac{\Gamma(x)\Gamma(y)}{\Gamma(x+y)}, ~~(Re(x), Re(y)>0).
\end{equation}
Let $\lambda \in (0,~1)$. Then, for $0<Re(s)<\frac{1-\lambda}{\lambda}$, we have
\begin{equation}\label{17}
\begin{split}
\Gamma_{\lambda}(s+1)=&\int_0 ^{\infty}e_{\lambda} ^{-t}t^sdt=\int_0 ^{\infty}(1+\lambda t)^{-\frac{1}{\lambda}}t^sdt\\
=&\left.(1+\lambda t)^{-\frac{1-\lambda}{\lambda}}t^s\frac{1}{\lambda-1}\right|_0 ^{\infty}+\frac{s}{1-\lambda}\int_0 ^{\infty}(1+\lambda t)^{-\frac{1-\lambda}{\lambda}}t^{s-1}dt\\
=&\frac{s}{1-\lambda}\int_0 ^{\infty} \left(1+\frac{\lambda}{1-\lambda}(1-\lambda)t\right)^{-\frac{1-\lambda}{\lambda}}((1-\lambda)t)^{s-1}\frac{1}{(1-\lambda)^{s-1}}dt\\
=&\frac{s}{(1-\lambda)^s}\int_0 ^{\infty} \left(1+\frac{\lambda}{1-\lambda}(1-
\lambda)t\right)^{-\frac{1-\lambda}{\lambda}}((1-\lambda)t)^{s-1}dt\\
=&\frac{s}{(1-\lambda)^{s+1}}\int_0 ^{\infty} \left(1+\frac{\lambda}{1-\lambda}y\right)^{-\frac{1-\lambda}{\lambda}}y^{s-1}dy\\
=&\frac{s}{(1-\lambda)^{s+1}}\Gamma_{\frac{\lambda}{1-\lambda}}(s).
\end{split}
\end{equation}
Therefore, by \eqref{17}, we obtain the following theorem.
\begin{theorem}\label{thm1}
Let $\lambda \in (0,~1)$. Then, for $0<Re(s)<\frac{1-\lambda}{\lambda}$, we have
\begin{equation*}
\Gamma_{\lambda}(s+1)=\frac{s}{(1-\lambda)^{s-1}}\Gamma_{\frac{\lambda}{1-\lambda}} (s).
\end{equation*}
\end{theorem}

Let $\lambda \in (0,~\frac{1}{2})$. Then, in view of  Theorem \ref{thm1}, for $1<Re(s)<\frac{1-\lambda}{\lambda}$, we get
\begin{equation}\label{18}
\begin{split}
\Gamma_{\lambda}(s+1)=&\frac{s}{(1-\lambda)^{s+1}}\Gamma_{\frac{\lambda}{1-\lambda}}(s)=\left(\frac{s}{(1-\lambda)^{s+1}}\right)\left(\frac{(1-\lambda)^s(s-1)}{(1-2\lambda)^s}\right)\Gamma_{\frac{\lambda}{1-2\lambda}}(s-1)\\
=&\frac{s(s-1)}{(1-\lambda)(1-2\lambda)^s}\Gamma_{\frac{\lambda}{1-2\lambda}}(s-1).
\end{split}
\end{equation}

Let $\lambda \in (0,~\frac{1}{3})$. Then, invoking Theorem 2.1 again, for $2<Re(s)<\frac{1-\lambda}{\lambda}$,  we have
\begin{equation}\label{19}
\begin{split}
\Gamma_{\lambda}(s+1)=&\frac{s(s-1)}{(1-\lambda)(1-2\lambda)^s}\Gamma_{\frac{\lambda}{1-2\lambda}}(s-1)\\
=&\frac{s(s-1)(s-2)(1-2\lambda)^{s-1}}{(1-\lambda)(1-2\lambda)^s(1-3\lambda)^{s-1}}\Gamma_{\frac{\lambda}{1-3\lambda}}(s-2)\\
=&\frac{s(s-1)(s-2)}{(1-\lambda)(1-2\lambda)(1-3\lambda)^{s-1}}\Gamma_{\frac{\lambda}{1-3\lambda}}(s-2).
\end{split}
\end{equation}
Let $\lambda \in (0,~\frac{1}{k+1})$. Continuing this process, for $k<Re(s)<\frac{1-\lambda}{\lambda}$,  we have
\begin{equation}\label{20}
\begin{split}
\Gamma_{\lambda}(s+1)=&\frac{s(s-1)\cdots(s-(k+1)+1)}{(1-\lambda)(1-2\lambda)\cdots(1-k\lambda)(1-(k+1)\lambda)^{s-k+1}}\Gamma_{\frac{\lambda}{1-(k+1)\lambda}}(s-(k+1)+1)\\
=&\frac{s(s-1)\cdots(s-(k+1)+1)}{(1-\lambda)(1-2\lambda)\cdots(1-k\lambda)(1-(k+1)\lambda)^{s-k+1}}\Gamma_{\frac{\lambda}{1-(k+1)\lambda}}(s-k).
\end{split}
\end{equation}

Therefore, by \eqref{20}, we obtain the following theorem.
\begin{theorem}\label{thm2}
Let $\lambda\in \left(0,~\frac{1}{k+1}\right)$, with $k\in{\mathbb{N}}$.  For $k<Re(s)<\frac{1-\lambda}{\lambda}$, we have
\begin{equation*}
\Gamma_{\lambda}(s+1)=\frac{s(s-1)\cdots(s-(k+1)+1)}{(1-\lambda)(1-2\lambda)\cdots(1-k\lambda)(1-(k+1)\lambda)^{s-k+1}}\Gamma_{\frac{\lambda}{1-(k+1)\lambda}}(s-k).
\end{equation*}
\end{theorem}

Note that
\begin{equation}\label{21}
\left(1-(k+1)\lambda\right)^{s-k+1}\frac{\Gamma_{\lambda}(s+1)}{\Gamma_{\frac{\lambda}{1-(k+1)\lambda}}(s-k)}=\frac{s(s-1)(s-2)\cdots(s-k)}{(1-\lambda)(1-2\lambda)\cdots(1-k\lambda)},
\end{equation}
where  $\lambda\in \left(0,~\frac{1}{k+1}\right)$, with $k\in{\mathbb{N}}$, and $k<Re(s)<\frac{1-\lambda}{\lambda}$.

Let us take $s=k+1$. Then, by Theorem \ref{thm2}, we get
\begin{equation}\label{22}
\Gamma_{\lambda}(k+2)=\frac{(k+1)!}{(1-\lambda)(1-2\lambda)\cdots(1-k\lambda)(1-(k+1)\lambda)^2}\Gamma_{\frac{\lambda}{1-(k+1)\lambda}}(1),
\end{equation}
where $\lambda \in (0,~\frac{1}{k+2})$.

Now, we observe that
\begin{equation}\label{23}
\begin{split}
\Gamma_{\frac{\lambda}{1-(k+1)\lambda}}(1)=&\int_0 ^{\infty}\left(1+\frac{\lambda}{1-(k+1)\lambda}t\right)^{-\frac{1-(k+1)\lambda}{\lambda}}dt\\
=&\left.\left(\frac{\lambda}{\lambda-1+(k+1)\lambda}\right)\left(\frac{1-(k+1)\lambda}{\lambda}\right)\left(1+\frac{\lambda t}{1-(k+1)\lambda}\right)^{-\frac{1-(k+2)\lambda}{\lambda}}\right|_0 ^{\infty}\\
=&\frac{1-(k+1)\lambda}{1-(k+2)\lambda},{\text{ for }}\lambda \in (0,~\frac{1}{k+2}).
\end{split}
\end{equation}
For $\lambda \in (0,~\frac{1}{k+2})$, by \eqref{20} and \eqref{21}, we get
\begin{equation}\label{24}
\Gamma_{\lambda}(k+2)=\frac{(k+1)!}{(1-\lambda)(1-2\lambda)\cdots(1-k\lambda)(1-(k+1)\lambda)(1-(k+2)\lambda)}.
\end{equation}
Therefore, by \eqref{24}, we obtain the following theorem.
\begin{theorem}\label{thm3}
For $k\in{\mathbb{N}}$ and $\lambda \in (0,~\frac{1}{k})$, we have
\begin{equation*}
\Gamma_{\lambda} (k)=\frac{(k-1)!}{(1-\lambda)(1-2\lambda)\cdots(1-k\lambda)}.
\end{equation*}
\end{theorem}

\section{Degenerate Laplace transform}

Let $\lambda \in (0,~\infty)$, and let $f(t)$ be a function defined for $t\geq0$. Then the integral
\begin{equation}\label{25}
\mathcal{L}_{\lambda}(f(t))=\int_0 ^{\infty} (1+\lambda t)^{-\frac{s}{\lambda}}f(t)dt
\end{equation}
is said to be the {\it{degenerate Laplace transform}} of $f$ if the integral converges, which is also denoted by $\mathcal{L}_{\lambda}(f(t))=F_{\lambda}(s)$.

From \eqref{25}, we note that
\begin{equation*}
\mathcal{L}_{\lambda}(\alpha f(t)+\beta g(t))=\alpha\mathcal{L}_{\lambda}(f(t))+\beta\mathcal{L}_{\lambda}(g(t)),
\end{equation*}
where $\alpha$ and $\beta$ are constant real numbers.

Now, we observe that
\begin{equation}\label{26}
\begin{split}
\mathcal{L}_{\lambda}(1)=&\int_0 ^{\infty}(1+\lambda t)^{-\frac{s}{\lambda}}dt=\lim_{b\rightarrow \infty}\int_0 ^b (1+\lambda t)^{-\frac{s}{\lambda}}dt\\
=&\frac{1}{s-\lambda}{\text{ if }}s>\lambda,
\end{split}
\end{equation}
and
\begin{equation}\label{27}
\begin{split}
\mathcal{L}_{\lambda}\left((1+\lambda t)^{-\frac{a}{\lambda}}\right)=&\int_0 ^{\infty}(1+\lambda t)^{-\frac{s+a}{\lambda}}dt\\
=&\frac{1}{s+a-\lambda},{\text{ if }}s>-a+\lambda.
\end{split}
\end{equation}
From \eqref{9} and \eqref{12}, we have
\begin{equation}\label{26}
\cos_{\lambda}(t)=\frac{e_{\lambda} ^{it}+e_{\lambda} ^{-it}}{2}=\frac{1}{2}\left((1+\lambda t)^{\frac{i}{\lambda}}+(1+\lambda t)^{-\frac{i}{\lambda}}\right),
\end{equation}
and
\begin{equation}\label{29}
\sin_{\lambda}(t)=\frac{e_{\lambda} ^{it}-e_{\lambda} ^{-it}}{2i}=\frac{1}{2i}\left((1+\lambda t)^{\frac{i}{\lambda}}-(1+\lambda t)^{-\frac{i}{\lambda}}\right).
\end{equation}
Note that
\begin{equation}\label{30}
\left(\cos_{\lambda}(t)\right)^{'}=\frac{d}{dt}\cos_{\lambda}(t)=-\frac{1}{1+\lambda t}\sin_{\lambda}(t),
\end{equation}
and
\begin{equation}\label{31}
\left(\sin_{\lambda}(t)\right)^{'}=\frac{d}{dt}\sin_{\lambda}(t)=\frac{1}{1+\lambda t}\cos_{\lambda}(t).
\end{equation}
As illustrations, we will compute the degenerate Laplace transforms of the degenerate trigonometric functions, degenerate hyperbolic functions and $t^n$. From now on, we will also refrain from stating any restrictions on $s$, with the understanding that $s$ is suitably restricted to guarantee the convergence of the relevant degenerate Laplace transfrom.
From \eqref{29}, we have
\begin{equation}\label{32}
\begin{split}
\mathcal{L}_{\lambda}\left(\cos_{\lambda}(at)\right)=&\int_0 ^{\infty}(1+\lambda t)^{-\frac{s}{\lambda}}\cos_{\lambda}(at)dt\\
=&\frac{1}{2}\int_0 ^{\infty}\left((1+\lambda t)^{-\frac{s-ai}{\lambda}}+(1+\lambda t)^{-\frac{s+ai}{\lambda}}\right)dt\\
=&\frac{1}{2}\left(\frac{1}{s-\lambda-ai}+\frac{1}{s-\lambda+ai}\right)=\frac{s-\lambda}{(s-\lambda)^2+a^2},
\end{split}
\end{equation}
and
\begin{equation}\label{33}
\begin{split}
\mathcal{L}_{\lambda}\left(\sin_{\lambda}(at)\right)
=&\frac{1}{2i}\int_0 ^{\infty}\left((1+\lambda)^{-\frac{s-ai}{\lambda}}-(1+\lambda t)^{-\frac{s+ai}{\lambda}}\right)dt\\
=&\frac{1}{2i}\left(\frac{1}{s-\lambda-ai}-\frac{1}{s-\lambda+ai}\right)\\
=&\frac{a}{(s-\lambda)^2+a^2}.
\end{split}
\end{equation}

Now, we define the degenerate hyperbolic cosine and degenerate hyperbolic sine functions as
\begin{equation}\label{34}
\cosh_{\lambda}(at)=\frac{1}{2}\left((1+\lambda t)^{\frac{a}{\lambda}}+(1+\lambda t)^{-\frac{a}{\lambda}}\right),
\end{equation}
and
\begin{equation}\label{35}
\sinh_{\lambda}(at)=\frac{1}{2}\left((1+\lambda t)^{\frac{a}{\lambda}}-(1+\lambda t)^{-\frac{a}{\lambda}}\right).
\end{equation}
Note that
\begin{equation}\label{36}
\left(\cosh_{\lambda}(t)\right)^{'}=\frac{d}{dt}\cosh_{\lambda}(t)=\frac{1}{1+\lambda t}\sinh_{\lambda}(t),
\end{equation}
and
\begin{equation}\label{37}
\left(\sinh_{\lambda}(t)\right)^{'}=\frac{d}{dt}\sinh_{\lambda}(t)=\frac{1}{1+\lambda t}\cosh_{\lambda}(t).
\end{equation}
By \eqref{34} and \eqref{35}, we get
\begin{equation}\label{38}
\begin{split}
\mathcal{L}_{\lambda}(\cosh_{\lambda}(at))=&\frac{1}{2}\int_0 ^{\infty}\left((1+\lambda t)^{-\frac{s-a}{\lambda}}+(1+\lambda t)^{-\frac{s+a}{\lambda}}\right)dt\\
=&\frac{1}{2}\left\{\frac{1}{s-a-\lambda}+\frac{1}{s+a-\lambda}\right\}\\
=&\frac{s-\lambda}{(s-\lambda)^2-a^2},
\end{split}
\end{equation}
and
\begin{equation}\label{39}
\begin{split}
\mathcal{L}_{\lambda}(\sinh_{\lambda}(at))=&\frac{1}{2}\int_0 ^{\infty}\left((1+\lambda t)^{-\frac{s-a}{\lambda}}-(1+\lambda t)^{-\frac{s+a}{\lambda}}\right)dt\\
=&\frac{1}{2}\left\{\frac{1}{s-a-\lambda}-\frac{1}{s+a-\lambda}\right\}\\
=&\frac{a}{(s-\lambda)^2-a^2}.
\end{split}
\end{equation}
Therefore, we obtain the following theorem
\begin{theorem}\textnormal{(Degenerate Laplace transforms of degenerate hyperbolic functions)}
\begin{equation*}
\mathcal{L}_{\lambda}(\cosh_{\lambda}(at))=\frac{s-\lambda}{(s-\lambda)^2-a^2},~
\mathcal{L}_{\lambda}(\sinh_{\lambda}(at))=\frac{a}{(s-\lambda)^2-a^2}.
\end{equation*}
\end{theorem}

Let $s>(n+1)\lambda$, with $n\in{\mathbb{N}}$. Then, as we saw in \eqref{15} and \eqref{16}, we have
\begin{equation}
\mathcal{L}_{\lambda}(t^n)=\int_0 ^{\infty}(1+\lambda t)^{-\frac{s}{\lambda}}t^ndt=\lambda^{-n-1}B(n+1,\frac{s}{\lambda}-n-1),
\end{equation}
and
\begin{equation}\label{40}
\begin{split}
&\mathcal{L}_{\lambda}(t^n)=\int_0 ^{\infty}(1+\lambda t)^{-\frac{s}{\lambda}}t^ndt\\
=&\left.-\frac{\lambda}{s-\lambda}\cdot\frac{1}{\lambda}(1+\lambda t)^{-\frac{s-\lambda}{\lambda}}t^n\right|_0 ^{\infty}+\frac{n}{s-\lambda}\int_0 ^{\infty}(1+\lambda t)^{-\frac{s-\lambda}{\lambda}}t^{n-1}dt\\
=&\frac{n}{s-\lambda}\int_0 ^{\infty} (1+\lambda t)^{-\frac{s-\lambda}{\lambda}}t^{n-1}dt.
\end{split}
\end{equation}
Let $s>(n+1)\lambda$. Then, by \eqref{40}, we get
\begin{equation}\label{41}
\begin{split}
\mathcal{L}_{\lambda}(t^n)=&\int_0 ^{\infty}(1+\lambda t)^{-\frac{s}{\lambda}}t^ndt=\frac{n}{s-\lambda}\int_0 ^{\infty} (1+\lambda t)^{-\frac{s-\lambda}{\lambda}}t^{n-1}dt\\
=&\frac{n(n-1)}{(s-\lambda)(s-2\lambda)}\int_0 ^{\infty}(1+\lambda t)^{-\frac{s-2\lambda}{\lambda}}t^{n-2}dt=\cdots\\
=&\frac{n!}{(s-\lambda)(s-2\lambda)\cdots(s-n\lambda)}\int_0 ^{\infty}(1+\lambda t)^{-\frac{s-n\lambda}{\lambda}}dt\\
=&\frac{n!}{(s-\lambda)(s-2\lambda)\cdots(s-n\lambda)(s-(n+1)\lambda)}\\
=&\frac{n!}{s^{n+1}}\frac{1}{\left(1-\frac{\lambda}{s}\right)\left(1-\frac{2\lambda}{s}\right)\cdots\left(1-\frac{(n+1)\lambda}{s}\right)}=\frac{1}{s^{n+1}}\Gamma_{\frac{\lambda}{s}}(n+1).
\end{split}
\end{equation}
Therefore, we obtain the following theorem.
\begin{theorem}
For $n\in{\mathbb{N}}$ and $s>(n+1)\lambda$, we have
\begin{equation*}
\begin{split}
\mathcal{L}_{\lambda}(t^n)=&\frac{n!}{(s-\lambda)(s-2\lambda)\cdots(s-n\lambda)(s-(n+1)\lambda)}\\
=&\frac{n!}{s^{n+1}}\frac{1}{\left(1-\frac{\lambda}{s}\right)\left(1-\frac{2\lambda}{s}\right)\cdots\left(1-\frac{(n+1)\lambda}{s}\right)}.
\end{split}
\end{equation*}
Moreover,
\begin{equation*}
{\mathcal{L}}_{\lambda}(t^n)=\frac{1}{s^{n+1}}\Gamma_{\frac{\lambda}{s}}(n+1).
\end{equation*}
\end{theorem}

Let $\alpha\in{\mathbb{R}}$ with $\alpha>-1$. Then, for $s>(\alpha+1)\lambda$, we have
\begin{equation}\label{42}
\begin{split}
\mathcal{L}_{\lambda}(t^{\alpha})=&\int_0 ^{\infty}(1+\lambda t)^{-\frac{s}{\lambda}}t^{\alpha}dt=\int_0 ^{\infty} \left(1+\frac{\lambda}{s}st\right)^{-\frac{s}{\lambda}}t^{\alpha}dt\\
=&\frac{1}{s^{\alpha}}\int_0 ^{\infty}\left(1+\frac{\lambda}{s}st\right)^{-\frac{s}{\lambda}}(st)^{\alpha}dt\\
=&\frac{1}{s^{\alpha+1}}\int_0 ^{\infty} \left(1+\frac{\lambda}{s}y\right)^{-\frac{s}{\lambda}}y^{\alpha}dy\\
=&\frac{1}{s^{\alpha+1}}\Gamma_{\frac{\lambda}{s}}(\alpha+1).
\end{split}
\end{equation}

The Laplace transforms are used in solving certain types of differential equations by reducing them to algebra problems. They are also used in such diverse areas as circular analysis, proportional-integral-derivative controllers, DC moter speed control systems and DC moter position control systems. It is expected that the degenerate Laplace transforms will find applications not only in mathematics but also in some applied areas. Next, as in the case of the usual Laplace transform, we would like to derive formulas for the degenerate Laplace transform of derivatives of functions.

We observe that
\begin{equation}\label{43}
\begin{split}
\mathcal{L}_{\lambda}\left(f^{(1)}(t)\right)=&\mathcal{L}_{\lambda}\left(\frac{d}{dt}f(t)\right)=\int_0 ^{\infty}(1+\lambda t)^{-\frac{s}{\lambda}}f^{'}(t)dt\\
=&-f(0)+s\mathcal{L}_{\lambda}\left((1+\lambda t)^{-1}f(t)\right).
\end{split}
\end{equation}
Note that
\begin{equation}\label{44}
\mathcal{L}_{\lambda}\left((1+\lambda t)^{-1}f^{'}(t)\right)=-f(0)+(s+\lambda)\mathcal{L}_{\lambda}\left((1+\lambda t)^{-2}f(t)\right).
\end{equation}
From \eqref{41} and \eqref{44}, we have
\begin{equation}\label{45}
\begin{split}
\mathcal{L}_{\lambda}\left(f^{(2)}(t)\right)=&\mathcal{L}_{\lambda}\left(\left(\frac{d}{dt}\right)^2f(t)\right)=-f^{'}(0)+s\mathcal{L}_{\lambda}\left((1+\lambda t)^{-1}f^{'}(t)\right)\\
=&-f^{'}(0)-sf(0)+s(s+\lambda)\mathcal{L}_{\lambda}\left((1+\lambda t)^{-2}f(t)\right).
\end{split}
\end{equation}
Continuing this process, we have
\begin{equation}\label{46}
\begin{split}
&\mathcal{L}_{\lambda}\left(f^{(n)}(t)\right)=\mathcal{L}_{\lambda}\left(\left(\frac{d}{dt}\right)^nf(t)\right)\\
=&s(s+\lambda)\cdots(s+(n-1)\lambda)\mathcal{L}_{\lambda}\left((1+\lambda t)^{-n}f(t)\right)-\sum_{i=0} ^{n-1}f^{(i)}(0)\left(\prod_{l=1} ^{n-i-1}s+(l-1)\lambda\right).
\end{split}
\end{equation}
A function $f(t)$ is said to be of {\it{degenerate exponential order $C$}} if there exist $C,M>0$ and $T>0$ such that
\begin{equation*}
\left|f(t)\right|\leq M(1+\lambda t)^{\frac{C}{\lambda}}{\text{ for all }}t>T.
\end{equation*}
If $f(t)$ is piecewise continuous on $(0,~\infty)$ and of degenerate exponential order $C$, then $\mathcal{L}_{\lambda}(f(t))$ exists for $s>C+\lambda$.

Therefore, we obtain the following theorem.
\begin{theorem}\label{thm6}
If $f,~f^{(1)},\cdots,f^{(n-1)}$ are continuous on $(0,~\infty)$ and are of degenerate exponential order and if $f^{(n)}(t)$ is piecewise continuous on $(0,~\infty)$, then
\begin{equation*}
\mathcal{L}_{\lambda}\left(f^{(n)}(t)\right)=s(s+\lambda)\cdots(s+(n-1)\lambda)\mathcal{L}_{\lambda}\left((1+\lambda t)^{-n}f(t)\right)-\sum_{i=0} ^{n-1}f^{(i)}(0)\left(\prod_{l=1} ^{n-i-1}s+(l-1)\lambda\right),
\end{equation*}
where $f^{(n)}(t)=\left(\frac{d}{dt}\right)^nf(t)$.
\end{theorem}

Let
\begin{equation}\label{47}
F_{\lambda}(s)=\int_0 ^{\infty}(1+\lambda t)^{-\frac{s}{\lambda}}f(t)dt.
\end{equation}
Then, by \eqref{47}, we get
\begin{equation}\label{48}
\begin{split}
\frac{d}{ds}F_{\lambda}(s)=&-\int_0 ^{\infty}\log(1+\lambda t)^{\frac{1}{\lambda}}(1+\lambda t)^{-\frac{s}{\lambda}}f(t)dt\\
=&-\lambda^{-1}\mathcal{L}_{\lambda}\left(\log(1+\lambda t)f(t)\right),
\end{split}
\end{equation}
and
\begin{equation}\label{49}
\begin{split}
\left(\frac{d}{ds}\right)^2F_{\lambda}(s)=&(-1)^2\int_0 ^{\infty}\left(\log(1+\lambda t)^{\frac{1}{\lambda}}\right)^2(1+\lambda t)^{-\frac{s}{\lambda}}f(t)dt\\
=&(-1)^2\lambda^{-2}\int_0 ^{\infty}\left(\log(1+\lambda t)\right)^2f(t)(1+\lambda t)^{-\frac{s}{\lambda}}dt\\
=&(-1)^2\lambda^{-2}\mathcal{L}_{\lambda}\left(\left(\log(1+\lambda t)\right)^2f(t)\right).
\end{split}
\end{equation}
Continuing this process, we have
\begin{equation}\label{50}
\begin{split}
\left(\frac{d}{ds}\right)^nF_{\lambda}(s)=&(-1)^n\int_0 ^{\infty}\left(\log(1+\lambda t)^{\frac{1}{\lambda}}\right)^n(1+\lambda t)^{-\frac{s}{\lambda}}f(t)dt\\
=&(-1)^n\lambda^{-n}\int_0 ^{\infty}\left(\log(1+\lambda t)\right)^nf(t)(1+\lambda t)^{-\frac{s}{\lambda}}dt\\
=&(-1)^n\lambda^{-n}\mathcal{L}_{\lambda}\left(\left(\log(1+\lambda t)\right)^nf(t)\right).
\end{split}
\end{equation}
Therefore, we obtain the following theorem.
\begin{theorem}\label{thm7}
For $n\in{\mathbb{N}}$, we have
\begin{equation*}
\mathcal{L}_{\lambda}\left(\left(\log(1+\lambda t)\right)^nf(t)\right)=(-1)^n\lambda^n\left(\frac{d}{ds}\right)^nF_{\lambda}(s).
\end{equation*}
\end{theorem}
By Theorem \ref{thm7}, we note that
\begin{equation}\label{51}
\begin{split}
\mathcal{L}_{\lambda}\left((1+\lambda t)^{\frac{a}{\lambda}}\left(\log(1+\lambda t)^n\right)\right)=&(-1)^n\lambda^n\left(\frac{d}{ds}\right)^n\left(\frac{1}{s-a-\lambda}\right)\\
=&\lambda^n n!\left(\frac{1}{s-a-\lambda}\right)^{n+1}.
\end{split}
\end{equation}
By Taylor expansion, we get
\begin{equation}\label{52}
\begin{split}
\mathcal{L}_{\lambda}\left(f(t)(1+\lambda t)^{\frac{a}{\lambda}}\right)=&\sum_{n=0} ^{\infty}\frac{a^n\lambda^{-n}}{n!}\mathcal{L}_{\lambda}\left(\left(\log(1+\lambda t)\right)^nf(t)\right)\\
=&\sum_{n=0} ^{\infty}\frac{(-a)^n}{n!}\left(\frac{d}{ds}\right)^nF_{\lambda}(s).
\end{split}
\end{equation}

\end{document}